\documentclass{article}

\usepackage{a4wide}
\usepackage{amsmath}
\usepackage{amssymb}
\usepackage{amsthm}
\usepackage{txfonts}
\usepackage{mathtools}
\usepackage{mdwlist}
\usepackage{hyperref}
\usepackage[svgnames]{xcolor}

\newtheorem{remark}{Remark}%\mumberwithin{remark}{section}
\newtheorem{lemma}{Lemma}%\mumberwithin{lemma}{section}
\newtheorem{theorem}{Theorem}%\mumberwithin{theorem}{section}
\newtheorem{ass}{Assumption}
\newtheorem*{alg*}{Algorithm}

\renewcommand{\eqref}[1]{\hyperref[#1]{(\ref{#1})}}

\newcommand{\bs}[1]{\boldsymbol{#1}}
\newcommand{\pt}{\partial}
\newcommand{\laplace}{\Delta}
\newcommand{\eps}{\varepsilon}
\newcommand{\RR}{\mathbb{R}}
\newcommand{\PP}{\Pi}
\renewcommand{\phi}{\varphi}
\renewcommand{\rho}{\varrho}
\renewcommand{\theta}{\vartheta}
\newcommand{\scal}[2]{\left(#1,#2\right)}
\newcommand{\dual}[2]{\left\langle#1,#2\right\rangle}
\newcommand{\bigscal}[2]{\bigl(#1,#2\bigr)}
\newcommand{\bigdual}[2]{\bigl\langle#1,#2\bigr\rangle}
\newcommand{\m}[1]{\mathcal{#1}}
\newcommand{\mL}{\mathcal{L}}

\newcommand{\NN}{\mathbb{N}}
\newcommand{\E}{\mathrm{e}}
\newcommand{\D}{\mathrm{  d}}
\newcommand{\abs}[1]{\left|#1\right|}
\newcommand{\norm}[1]{\left\|#1\right\|}
\newcommand{\bignorm}[1]{\big\|#1\big\|}

\newif\ifNUMERICS\NUMERICSfalse

\hypersetup{
  pdftitle    = {Maximum-norm a posteriori error bounds for parabolic equations
       discretised by the extrapolated Euler method in time and FEM in space},
  pdfsubject  = {Numerical Analysis},
  pdfauthor   = {Torsten Lin\ss, Goran Radojev},
  pdfkeywords = {maximum-norm a posteriori error estimates, extrapolation,
  upwind sheme, singular perturbation, convection-diffusion, adaptivity},
  pdfcreator  = {pdflatex},
  pdfproducer = {LaTeX with hyperref},
    colorlinks,
    linkcolor={red!90!black},
    citecolor={blue!90!black},
    urlcolor={green!90!black}
}

\title{Maximum-norm a posteriori error bounds for parabolic equations
       discretised by the extrapolated Euler method in time and FEM in space}

\author{Torsten Lin\ss\thanks{Fakult\"at f\"ur Mathematik und Informatik,
        FernUniversit\"at in Hagen,
        Universit\"atsstra{\ss}e 11,
        58095 Hagen,
        Germany,
        \texttt{torsten.linss@fernuni-hagen.de}}
   \and Goran Radojev\thanks{Department of Mathematics and Computer Science, Faculty of Sciences,
        University of Novi Sad, Trg Dositeja Obradovi\'ca~4, 21000 Novi Sad,
        Serbia,
        \texttt{goran.radojev@dmi.uns.ac.rs}.
        }}

\begin{document}

\maketitle

\begin{abstract}
  A class of linear parabolic equations is considered.
  We derive a framework for the a posteriori error analysis
  of time discretisations by Richardson extrapolation of arbitrary order
  combined with finite element discretisations in space.
  We use the idea of elliptic reconstructions and certain bounds
  for the Green’s function of the parabolic operator.
  The crucial point in the analysis is the design of suitable polynomial
  reconstructions in time from approximations that are given only in
  the mesh points.
  
  \emph{Keywords:}
  parabolic problems, maximum-norm a posteriori error estimates,
  backward Euler, extrapolation, elliptic reconstructions, Green's function.

  \emph{AMS subject classification (2000):} 65M15, 65M50, 65M60.
\end{abstract}

%%%%%%%%%%%%%%%%%%%%%%%%%%%%%%%%%%%%%%%%%%%%%%%%%%%%%%%%%%%%%%%%%%%%%%%%%%%%%%%%
\section{Introduction}

Given a second-order linear elliptic operator $\m{L}$ in an open spatial
domain $\Omega\subset\RR^n$ with Lipschitz boundary, we consider the linear
parabolic equation
\begin{subequations}\label{problem}
\begin{alignat}{2}
 \m{K}u \coloneqq \pt_t u  + \m{L} u & = f\,,
            &\quad& \text{in} \quad Q\coloneqq\Omega \times (0,T],\\
 \intertext{subject to the initial condition}
   u(x,0) & = u^0(x)\,, && \text{for} \quad x\in\bar{\Omega}, \\
 \intertext{and the homogeneous Dirichlet boundary condition}
   u(x,t) & = 0\,, && \text{for} \quad (x,t)\in \pt\Omega \times [0,T].
\end{alignat}
\end{subequations}
%
%Precise assumptions on the data will be given later.

There is a vast literature dealing with numerical methods for approximating
parabolic problems.
One classical example is Thom\'ee's monograph \cite{MR2249024} which gives a
comprehensive overview of the FEM for~\eqref{problem} and related equations.
The subject is also elaborated on in various textbooks, e.g.~\cite{MR2362757,
MR3015004}.
The majority of these publications focus on \emph{a priori} error estimation,
i.e. convergence results given in terms of the mesh size and the regularity of
the \emph{exact} solution.

More recently, the derivation of \emph{a posteriori} error bounds has attracted
the attention of many researchers.
This approach yields error bounds in terms of the \emph{computed} numerical
solution, and therefore yields computable bounds on the error.

Most publications to date study the error in $L_2$-norms or in energy norms
naturally induced by the variational formulation of the problem, see
e.g.~\cite{MR695602,MR695603,MR1043607}.
In contrast, our focus is on a posteriori error bounds in the \emph{maximum norm}.
The first such results were probably given in~\cite{MR1335652}, however the
proofs were deferred to a later paper which seems to have never been published.

A key publication is~\cite{MR2034895} by Makridakis and Nochetto who introduced
the concept of \emph{elliptic reconstructions}.
This idea was used in a number of publications to study the FEM combined with
various time discretisations: backward Euler, Crank Nicolson and discontinuous
Galerkin~\cite{MR2519598, MR2629992, MR3056758}.
The techniques in those papers are all tailored to the particular method(s)
under consideration and therefore differ to some extend.
Moreover, different --- although related --- stability results have been used.

In~\cite{MR4755064} a common framework for the maximum-norm a posteriori
error analysis of second order time discretisations was derived and
applied to various time discretisations: backward Euler, Crank-Nicolson,
a two-step extrapolation method based on backward Euler, the
BDF-2 method, the Lobatto IIIC method and a two-stage SDIRK method.

In the present paper we shall study the use of Richardson extrapolations of
arbitrary order \mbox{$L\ge2$} for the time discretisation, while in space
any finite-element method might be used.
The case \mbox{$L=2$} was studied in~\cite{MR4755064}.
There the point-wise approximations in the nodes of the temporal mesh were
joined by ``straight lines'' (i.e. $\Pi_1$-interpolation) to obtain a
continuous approximation on the whole time interval \mbox{$[0,T]$}.
For \mbox{$L>2$} the design of this \emph{temporal reconstruction} is more
delicate.
To match to order $L$ of the method, on each subinterval the reconstruction
should use a polynomial of degree $L-1$.
However, only in the mesh points the approximation is of order $L$, all
other auxilliary approximations constructed during the extrapolation process
are of lower order and therefore useless for our purpose.
This crucial issue will be addressed in~\autoref{sssect:apost-para-temp_rec}.

The paper is organised as follows.
In \autoref{sect:weak} we specify our assumptions on the data of
problem~\eqref{problem} and recapitulate certain aspects of the existence theory.
\autoref{sect:discr} contains details of the time discretisation by
the Richardson extrapolation technique applied to the implicit Euler method.

In \autoref{sect:ana} we develop our a posteriori error analysis.
We formulate the assumptions for the analysis in \autoref{ssect:apost-ell}
(existence of error estimators for the elliptic problems) and
in \autoref{ssect:green} (bounds on the Green's function of the parabolic problem).
In \autoref{ssect:reco} the concept of elliptic reconstructions is recalled,
while the main result  is derived in \autoref{ssect:apost-para}.
\ifNUMERICS
{\color{blue}
Finally, numerical results are presented in \autoref{sect:numer} to illustrate our
theoretical findings.
}
\fi

\section{Weak formulation}
\label{sect:weak}

We shall study~\eqref{problem} in its standard variational form,
cf.~\cite[\S5.1.1]{MR2362757}.
The appropriate Gelfand triple consists of the spaces
\begin{gather*}
  V = H_0^1(\Omega), \quad H= L_2(\Omega) \quad\text{and}\quad
  V^* = H^{-1}(\Omega)\,.
\end{gather*}
Moreover, by \mbox{$a(\cdot,\cdot) \colon V \times V \to \RR$}
we denote the bilinear form associated with the elliptic operator $\m{L}$,
while $\dual{\cdot}{\cdot}\colon V^*\times V \to \RR$ is the duality
pairing and $\scal{\cdot}{\cdot}\colon H\times H\to\RR$ is the scalar
product in~$H$.
Also we shall denote by $\norm{\cdot}_{q,D}$ the standard
norm in $L_q(D)$, \mbox{$q\in[1,\infty]$}.

The solution $u$ of~\eqref{problem} may be considered as a mapping
\mbox{$[0,T] \to V\colon t \mapsto u(t)$}, and we will denote its
(temporal) derivative by $u'$ (and $\pt_t u$).
Let
\begin{gather*}
  W_2^1(0,T;V,H) \coloneqq \left\{ v \in L_2(0,T;V) \colon
         v'\in L_2(0,T;V^*)\right\}\,.
\end{gather*}
Then our variational formulation of~\eqref{problem} reads:
Given \mbox{$u^0\in H$} and \mbox{$F\in L_2(0,T;V^*)$},
find \mbox{$u\in W_2^1(0,T;V,H)$} such that
\begin{subequations}\label{weak}
\begin{alignat}{2}
   \frac{\D}{\D t} \bigscal{u(t)}{\chi} + a\bigscal{u(t)}{\chi}
            & = \bigdual{F(t)}{\chi} \quad \forall \chi\in V,
      \ \ t\in(0,T], \\
  \intertext{and}
    u(0)=u^0.
\end{alignat}
\end{subequations}
This problem possesses a unique solution.
In the sequell we shall assume that the source term $F$ is more
regular and can be represented as \mbox{$\dual{F(t)}{\chi}=\scal{f}{\chi}$},
\mbox{$\forall v\in V$}, with a function \mbox{$f\in L_2(0,T;H)$}

We are interested in maximum-norm error estimates.
Therefore, we have to make further
assumptions on the data to ensure that the solution can be evaluated pointwise.
To this end, we assume that the intial and boundary data satisfy the zero-th
order compatibility condition, i.e. $u^0 = 0$ on $\pt\Omega$, and that $u^0$
is Hölder continuous in $\bar\Omega$.
Under standard assumptions on $f$ and $\m{L}$, problem~\eqref{problem} possesses
a unique solution that is continuous on $\bar{Q}$; see~\cite[\S5, Theorem 6.4]{0174.15403}.

\section{Time discretisation}
\label{sect:discr}

Let a mesh in time be given by
\begin{gather*}
  \omega \colon 0=t_0<t_1<\ldots<t_M=T,
\end{gather*}
with mesh intervals \mbox{$I_j\coloneqq(t_{j-1},t_j)$} and
step sizes \mbox{$\tau_j\coloneqq t_j-t_{j-1}$}, \mbox{$j=1,2,\dots,M$}.
Furthermore, for \mbox{$\xi\in[0,1]$},
set \mbox{$t_{j-\xi}\coloneqq t_j-\xi\tau_j$}, and
for any function \mbox{$v\in C^0\left([0,T],L_2(\Omega)\right)$}
let \mbox{$v^{j-\xi} \coloneqq v(t_{j-\xi})$}, \mbox{$j=0,1,\dots,M$}.

%Let $V_h$ be a finite dimentional (FE-)subspace of $V$.
%Let us be given a time stepping procedure that generates a sequence
%\mbox{$\bigl\{u_h^j\bigr\}_{j=0}^M$} of approximations
%to $u$ on the mesh $\omega$, i.e., \mbox{$u(t_j)\approx u_h^j \in V_h$},
%\mbox{$j=0,\dots,M$}.
%Using any means of interpolation this sequence can be extended to a function
%\mbox{$u_h\in C^0\bigl([0,T],V_h\bigr)$} such that \mbox{$u_h(t_j)=u_h^j$},
%\mbox{$j=0,\dots,M$}.
%We shall refer to this function as the \emph{temporal reconstruction} of
%the $u_h^j$.
%Its construction will later be taylored to the particular time stepping method
%under consideration.
%
%\begin{gather*}
%  v^0 \coloneqq v(\cdot,0) \quad\text{and}\quad
%  v^{j-\mu} \coloneqq v(\cdot,t_j-\mu \tau_j)\, \ \ j=1,2,\dots,M, \ \
%    \mu\in[0,1]
%\end{gather*}
%\footnote{
%We also set
%%
%\begin{gather*}
%  \delta_t v^j \coloneqq \frac{v^j-v^{j-1}}{\tau_j}\,, \ \ j=1,2,\dots,M.
%\end{gather*}
%}

% We shall consider an arbitrary time-stepping procedure that generates as
% sequence \mbox{$u^j_h\in V_h$}, \mbox{$j=0,1,\dots,M$}, of approximations
% to the solution~$u$ of~\eqref{problem} at time $t_j$,
% i.e. \mbox{$u^j_h\approx u(t_j)$}.

For the temporal discretisation of~\eqref{weak} we shall use Richardson
extrapolation with $L\ge2$ extrapolation steps applied to the
backward Euler method.
We slightly modify~\eqref{weak} by replacing $f$ by $\hat{f}$, where
$\hat{f}\vert_{I_j}\in \PP_{L-1}$, $j=1,\dots,M$, interpolates $f$ at
the mesh points and on \mbox{$L-2$} further points on each mesh interval
$I_j$.

For the spatial discretisation,
let $V_h$ be a finite dimentional (FE-)subspace of $V$, and
let $a_h\scal{\cdot}{\cdot}$ and $\scal{\cdot}{\cdot}_h$
be approximations of the bilinear form $a\scal{\cdot}{\cdot}$
and of the scalar product $\scal{\cdot}{\cdot}$.
These may involve quadrature, for example.
Furthermore, let $u_h^0\in V_h$ be an approximation of the initial
condition $u^0$,

Then our discretisation is as follows:
\begin{alg*}\mbox{}

\noindent
\textbf{S1: Initialise constants} \\[1.5ex]
\mbox{}\qquad
  $\displaystyle
      \alpha_\ell \coloneqq \frac{\abs{\bs{\mathrm{e}}^\ell \vert \bs{B}}}
                                 {\abs{\bs{1}^{\phantom{\ell}} \vert \bs{B}}}\,, \ \
      \text{where} \ \
      \bs{B}\coloneqq \left(b_{ik}\right) \in \RR^{L,L-1}\,, \ \ b_{im}=1/i^m.
  $\\[1.5ex]
\textbf{S2: Define initial approximations} \\[1ex]
\textsc{for} $\ell=1,2,\dots,L$ \textsc{do} \\[1.5ex]
  \mbox{}\qquad
  $ v_\ell^0 = u_h^0$ \\[1.5ex]
\textsc{endfor} // $\ell$\\

\noindent
\textbf{S3: Time stepping procedure} \\[1.5ex]
\textsc{for} $j=1,2,\dots,M$ \textsc{do} \\[1.5ex]
  \mbox{}\qquad
  \textsc{for} $\ell=1,2,\dots,L$ \textsc{do} \\[1.5ex]
    \mbox{}\qquad\qquad
    \textsc{for} $k=\ell-1,\dots,0$ \textsc{do} \\[1.5ex]
      \mbox{}\qquad\qquad\qquad
      Find $v_\ell^{j-k/\ell}$ such that \\[1.5ex]
      \mbox{}\qquad\qquad\qquad
      $\displaystyle
          \scal{\frac{v_\ell^{j-k/\ell} - v_\ell^{j-(k+1)/\ell}}{\tau_j/\ell}}{\chi}_h
               + a_h\scal{v_\ell^{j-k/\ell}}{\chi} = \scal{\hat{f}^{j-k/\ell}}{\chi}_h\,, \ \ \forall \ \chi \in V_h
      $ \\
    \mbox{}\qquad\qquad
    \textsc{endfor} // $k$ \\[1.5ex]
  \mbox{}\qquad
  \textsc{endfor} // $\ell$ \\[1.5ex]
  \mbox{}\qquad
    $\displaystyle
       u_h^j \coloneqq \sum_{\ell=1}^L \alpha_\ell v_\ell^j
       \quad\quad \biggl(\approx u(t_j)\biggr)
    $ \\
\textsc{endfor} // $j$
\end{alg*}

\pagebreak

\section{A posteriori error analysis}
\label{sect:ana}

Our analysis uses three main ingredients:
\begin{itemize*}
  \item a posteriori error bounds for the elliptic problem $\m{L}y=g$,
    see \autoref{ssect:apost-ell},
  \item bounds for the \textsc{Green}'s function associated with the
    parabolic operator $\m{K}$, see \autoref{ssect:green} and
  \item the idea of elliptic reconstructions introduced by
    Makridakis and Nochetto~\cite{MR2034895}, see \autoref{ssect:reco}.
\end{itemize*}
In~\autoref{ssect:apost-para}, after these concepts have been reviewed,
we derive an a posteriori error estimator for the extrapolated Euler method
(and FEM in space) applied to~\eqref{problem}.
As mentioned earlier, the crucial point will be the design of the temporal
reconstruction of $u_h$.

\subsection{A posteriori error estimation for the elliptic problem}
\label{ssect:apost-ell}

Given $g\in H$, consider the elliptic boundary-value problem of
finding $y\in V$ such that
\begin{gather}\label{prob:ell}
   a\scal{y}{\chi} = \scal{g}{\chi}\,, \ \ \forall \chi\in V,
\end{gather}
and its discretisation of finding $y_h\in V_h$ such that
\begin{gather}\label{FEM:ell}
   a_h\scal{y_h}{\chi} = \scal{g}{\chi}_h\,, \ \ \forall \chi\in V_h.
\end{gather}

\begin{ass}\label{ass:ee}
  There exists an a posteriori error estimator $\eta$ for the
  FEM~\eqref{FEM:ell} applied to the elliptic problem~\eqref{prob:ell} with
  \begin{gather*}
    \norm{y_h-y}_{\infty,\Omega} \le \eta\bigl(y_h, g\bigr).
  \end{gather*}
\end{ass}

A few error estimators of this type are available in the literature.
We mention some of them.
\begin{itemize*}
  \item Nochetto et al.~\cite{05068766} study the semilinear problem
        $-\laplace u + g(\cdot,u)=0$ in up to three space dimensions.
        They give a posteriori error bounds for arbitrary order FEM
        on quasiuniform triangulations.
  \item Demlow and Kopteva~\cite{MR3520007} too consider arbitrary order FEM
        on quasiuniform triangulations, but for the singularly perturbed
        equation $-\eps^2\laplace u + g(\cdot,u)=0$.
        A posteriori error estimates are established that are robust in the
        perturbation parameter.
        Furthermore, in~\cite{MR3419888} for the same problem $P_1$-FEM on
        \emph{anisotropic} meshes are investigated, while
        Kopteva and Rankin~\cite{MR4625899} derive related results for
        the discontinuous Galerkin FEM.
  \item Convection dominated problems have been studied by Demlow, Franz
        and Kopteva~\cite{MR4648515}. A posteriori error bounds in the maximum
        norm are given for classical conforming and for stabilised FEM.
  \item In~\cite{MR2334045,MR3232628} arbitrary order FEM for
        the linear problem \mbox{$-\eps^2 u'' + ru = g$} in \mbox{$(0,1)$},
        \mbox{$u(0)=u(1)=0$} are considered.
        In contrast to the afore mentioned contributions all constants
        appearing in the error estimator are given explicitly.
\end{itemize*}

\subsection{Green's functions}
\label{ssect:green}

Let the \textsc{Green}'s function associated with $\m{K}$ and an arbitrary
point \mbox{$x\in\Omega$} be denoted by $\m{G}$,
Then for all \mbox{$\phi\in W^1_2\left(0,T;V,H\right)$}
\begin{gather}\label{green-rep}
  \phi(x,t) = \bigscal{\phi(0)}{\m{G}(t)}
              + \int_0^t \bigdual{\bigl(\m{K}\phi\bigr)(s)}{\m{G}(t-s)} \D s,
    \ \ \text{see~\cite{0144.34903}.}
\end{gather}
The \textsc{Green}'s function
\mbox{$\m{G} \colon \bar{\Omega}\times[0,T]\to\RR$}, \mbox{$t \in (0,T]$}.
solves, for fixed $x$,
\begin{gather*}
  \pt_t\m{G}+\m{L}^*\m{G} = 0, \ \ \text{in} \ \Omega\times\RR^+, \ \
  \m{G}\bigr|_{\pt\Omega} = 0, \ \ \m{G}(0) = \delta_x=\delta(\cdot-x)\,.
\end{gather*}

\begin{ass}\label{ass:green}
  There exist non-negative constants
  $\kappa_0$, $\kappa_1$, $\kappa_1'$ and $\gamma$ such that
  \begin{gather}\label{source:ass}
    \norm{\m{G}(t)}_{1,\Omega} \le \kappa_0\,\E^{-\gamma t} \eqqcolon \phi_0(t),
       \quad
    \norm{\pt_t \m{G}(t)}_{1,\Omega}
       \le \left(\frac{\kappa_1}{t} +\kappa_1'\right)
                  \,\E^{-\gamma t} \eqqcolon \phi_1(t),
  \end{gather}
  for all $x\in\bar\Omega$, $t\in[0,T]$.
\end{ass}

\ifNUMERICS
In \autoref{sect:numer} we will present numerical results for a
test problem that satisfies these assumptions.
\fi
A more detailed discussion of problem classes for which such results are
available is given in~\cite[\S2]{MR3720388},
see also~\cite[Appendix A]{MR3056758}.

\subsection{Elliptic reconstruction}
\label{ssect:reco}

Given approximations $u_h\in C^0\bigl([0,T],V_h\bigr)$ of $u$
and $\hat{f}\in C^0\bigl([0,T],L_2(\Omega)\bigr)$ of $f$, we define
$\psi\in C^0\bigl([0,T],V_h\bigr)$ by
\begin{gather}\label{reconstr-psi}
  \scal{\psi(t)}{\chi}_h = a_h\scal{u_h(t)}{\chi}
              - \scal{\hat{f}(t)}{\chi}_h \quad \forall \ \chi\in V_h\,,
      \ \ t\in[0,T].
\end{gather}
This can be reformulated as an ``elliptic'' problem:
\begin{gather}\label{elliptic:discr}
  a_h\scal{u_h(t)}{\chi} = \scal{\hat{f} + \psi(t)}{\chi}_h \quad
      \forall \ \chi\in V_h, \ \ t\in[0,T].
\end{gather}
Next, let $R(t)\in C^0\bigl([0,T],H_0^1(\Omega)\bigr)$ solve
\begin{gather}\label{reconstr}
  a\scal{R(t)}{\chi} = \scal{\hat{f}(t)+\psi(t)}{\chi}
    \quad \forall \ \chi\in H_0^1(\Omega)\,, \ \ t\in[0,T],
\end{gather}
or for short:
\begin{gather}\label{reco-LRfpsi}
  \m{L}R(t) = \hat{f}(t)+\psi(t)\,, \ \ t\in[0,T].
\end{gather}
The function $R$ is referred to as the elliptic reconstruction of
$u_h$, \cite{MR2034895}.

Now, $u_h(t)$ can be regarded as the finite-element approximation of $R(t)$
obtained by~\eqref{elliptic:discr}.
Its error \mbox{$E(t)\coloneqq R(t)-u_h(t)$}, \mbox{$t\in[0,T]$}, can be
bounded using the estimator from~\autoref{ssect:apost-ell}:
\begin{gather}\label{reconstr-esti}
  \norm{E(t)}_\infty =
  \norm{\bigl(u_h - R\bigr)(t)}_\infty
    \le \eta_\mathrm{ell}(t)
    \coloneqq \eta\left(u_h(t), \bigl(\hat{f}+\psi\bigr)(t)\right)\,,
      \ \ t\in[0,T].
\end{gather}

\begin{remark}
  By~\eqref{reconstr-psi} and the inner loop of the time stepping procedure
  we have for the Richardson extrapolation
  \begin{gather}\label{psi-comp}
     \psi(t_j) = \psi^j = - \sum_{\ell=1}^L \alpha_\ell
                  \frac{v_\ell^{j} - v_\ell^{j-1/\ell}}{\tau_j/\ell}
     \,, \ \ j=1,\dots,M.
  \end{gather}
\end{remark}

\subsection{A posteriori error estimation for the parabolic problem}
\label{ssect:apost-para}

We are now in a position to derive the main result of the paper.

\subsubsection{General framework}
\label{sssect:apost-para-gen}

We use the \textsc{Green}'s function representation~\eqref{green-rep} with
$\phi$ replaced by \mbox{$u-u_h$} to express the error at final
time $T$ and for any $x\in\Omega$ as
\begin{gather}\label{rep-error}
  \left(u(T)-u_h^M\right)(x) =
  \bigl(u-u_h\bigr)(x,T) = \bigscal{u^0-u_h^0}{\m{G}(T)}
              + \int_0^T \bigdual{\bigl(\m{K}(u-u_h\bigr)(s)}{\m{G}(T-s)} \D s.
\end{gather}
This involves the residual of $u_h$ in the differential equation
for which we have the representation
\begin{align*}
  \m{K}(u-u_h)(s)
    & = f(s) - \mL u_h(s) - \pt_t u_h(s) \\
    & = \bigl(f-\hat{f}\bigr)(s)
             % + \mL\left(R - u_h\right)(s)
             + \mL E(s)
         - \psi(t) - \pt_t u_h(s)\,, \ \ s\in [0,T] \setminus \omega_t,
\end{align*}
by~\eqref{reco-LRfpsi} and \mbox{$E=R-u_h$}.
This is substituted into~\eqref{rep-error} to give
\begin{gather}\label{error-rep-discr}
  \begin{split}
  \left(u(T)-u_h^M\right)(x)
     & = \scal{u^0-u^0_h}{\m{G}(T)}
          + \int_{0}^{T}\scal{\bigl(f-\hat{f}\bigr)(s)}{\m{G}(T-s)} \D s \\
     & \qquad
       + \int_0^T \dual{\mL E(s)}{\m{G}(T-s)} \D s
       - \int_0^T \dual{\left(\psi+\pt_t u_h\right)(s)}{\m{G}(T-s)} \D s\,.
  \end{split}
\end{gather}
The first two terms can be bounded using H\"older's inequality and
\eqref{source:ass}:
\begin{align}\label{est-init}
  \abs{\scal{u^0-u_h^0}{\m{G}(T)}}
     & \le \kappa_0 \sigma_0 \norm{u^0-u_h^0}_{\infty,\Omega}
  \intertext{and}\label{est-F}
  \abs{\int_{I_j} \scal{\bigl(f-\hat{f}\bigr)(s)}{\m{G}(T-s)}\D s}
     & \le \kappa_0 \sigma_j
              \int_{I_j} \norm{\bigl(f-\hat{f}\bigr)(s)}_{\infty,\Omega} \D s\,,
\end{align}
where we have set
\begin{gather}\label{def:sigma_j}
  \sigma_j \coloneqq \exp\left(-\gamma\left(T-t_j\right)\right)\,, \ \ j=0,\dots,M.
\end{gather}

The remaining two terms in~\eqref{error-rep-discr}, namely
\begin{gather*}
  \int_0^T \dual{\mL E(s)}{\m{G}(T-s)} \D s
  \quad \text{and} \quad
  \int_0^T \dual{\left(\psi+\pt_t u_h\right)(s)}{\m{G}(T-s)} \D s
\end{gather*}
have to be treated individually depending on the (piecewise) polynomial degree
of the temporal reconstruction $u_h$.

\subsubsection{Temporal reconstructions}
\label{sssect:apost-para-temp_rec}

The temporal reconstruction of the $u_h^j$, \mbox{$j=0,1,\dots,N$},
will be designed in a systematic manner.
Starting from a continuous piecewise $\PP_1$-interpolation,
on each interval $I_j$ we add polynomial bubbles.
A suitable basis for these are the integrated Legendre polynomials.

Let $P_i$ denote the Legendre polynomial of order \mbox{$i\in\NN_0$}.
Then the integrated Legendre polynomials $N_i$ are defined by
\begin{gather*}
  N_{i+1}(\xi) \coloneqq \int_{-1}^\xi P_i(\zeta) \D \zeta
    = \left(\xi^2-1\right) \frac{P_i'(\xi)}{i(i+1)}\,, \ \ \ i\in\NN.
\end{gather*}
Note that
\begin{gather*}
  N_{i+1} = \frac{P_{i+1}-P_{i-1}}{2i+1}\quad
  \text{and}\quad
  N_{i+1}(-1) = N_{i+1}(1) = 0\,, \ \ \ i\in\NN.
\end{gather*}

Next, we define the temporal reconstruction for the $L$-step Richardson
extrapolation method, \mbox{$L\ge2$}.
To this end, let both the reconstruction $u_h$ and the interpolant $\hat{f}$
of $f$ be of polynomial degree $L-1$ on each~$\bar{I}_j$.
This degree matches the order of the time discretisation.
Because of~\eqref{reconstr-psi}, both the function $\psi$ and the elliptic
reconstruction $R$ are of polynomial degree $L-1$ on $\bar{I}_j$ too.

Any function \mbox{$\phi\in C^0\left([0,T],V\right)$} that is of piecewise
polynomial degree $L-1$ can be written as
\begin{gather}\label{t-poly-p}
  \phi = \Delta^0 \phi^j + \Delta^1 \phi^j P_1\circ \xi_j
              + \sum_{i=2}^{L-1} \Delta^i\phi^j N_{i}\circ \xi_j
     \ \ \ \text{in} \ \bar{I}_j\,, \ \ j=1,\dots,M.
\end{gather}
where
\begin{gather*}
    \Delta^0 \coloneqq \frac{\phi^j+\phi^{j-1}}{2}\,,\qquad
    \Delta^1 \phi^j \coloneqq \frac{\phi^j-\phi^{j-1}}{2}\,,
\end{gather*}
while the \mbox{$\xi_j\colon \bar{I_j} \to [-1,1]$} are affine maps from the time
intervals~$\bar{I}_j$ onto the reference interval $[-1,1]$:
\begin{gather*}
  \xi_j\colon \bar{I_j} \to [-1,1]\colon
     s \mapsto \xi_j(s) \coloneqq \frac{2\bigl(s-t_{j-1/2}\bigr)}{\tau_j}\,,
     \ \ \ j=1,\dots,M.
\end{gather*}
Note that the $\Delta^i\phi^j$, \mbox{$i=2,\dots,p$}, uniquely define $\phi$
on $\bar{I}_j$.
On the other hand, prescribing the values of $\phi$ at additional
\mbox{$L-2$} points in $I_j$ will determine the $\Delta^i\phi^j$.
% In the latter case the $\Delta^i\phi^j$ may be interpreted as difference
% quotients of order $i$.

%Adapting the notation from~\eqref{t-poly-p}, we have
%%
%\begin{gather*}
%  u_h(t) = \Delta^0 u_h^j + \Delta^1 u_h^j P_1\left(\xi_j(t)\right)
%              + \sum_{i=2}^{L-1} \Delta^i u_h^j N_{i}\left(\xi_j(t)\right)\,,
%    \ \ \ t\in\bar{I}_j, \ \ j=1,\dots,M,
% \intertext{and}
%  \hat{f}(t) = \Delta^0 \hat{f}^j + \Delta^1 \hat{f}^j P_1\left(\xi_j(t)\right)
%              + \sum_{i=2}^{L-1} \Delta^i \hat{f}^j N_{i}\left(\xi_j(t)\right)\,,
%    \ \ \ t\in\bar{I}_j, \ \ j=1,\dots,M.
%\end{gather*}

While the $\Delta^i\hat{f}^j$ are uniquely defined by the particular
interpolation of $f$ used in the implementation of the Richardson extrapolation
procedure,
the $\Delta^i u_h^j$, $i=2,\dots,L-1$, have to be chosen carefully to obtain
suitable a posteriori error bounds.

As stated, on each $I_j$ the functions/reconstructions $\hat{f}$, $u_h$,
$\psi$ and $R$ are elements of $\PP_{L-1}$, and
in view of~\eqref{reconstr-psi}, \eqref{elliptic:discr} and~\eqref{reconstr},
we have
\begin{gather}\label{reconstr-delta-psi}
 a_h\scal{\Delta^iu_h^j}{\chi} = 
  \scal{\Delta^i\psi^j + \Delta^i \hat{f}^j}{\chi}_h
    \quad \forall \ \chi\in H_0^1(\Omega)\,, \ \ i=0,1,\dots,L-1,
\end{gather}
and
\begin{gather*}
  a\scal{\Delta^i R^j}{\chi} = \scal{\Delta^i \hat{f}^j+\Delta^i \psi^j}{\chi}
    \quad \forall \ \chi\in H_0^1(\Omega)\,, \ \ i=0,1,\dots,L-1.
\end{gather*}
The latter and~\eqref{reconstr-esti} imply
\begin{gather}\label{reconstr-esti-delta}
  \norm{\Delta^i E^j}_\infty
  \le \eta^j_\mathrm{ell,i}
  \coloneqq \eta\left(\Delta^i u_h^j, \Delta^i \hat{f}^j+\Delta^i\psi^j\right)\,,
    \quad i=0,1,\dots,L-1, \ \ j=1,\dots,M.
\end{gather}

We now define the \mbox{$\Delta^i u_h^j$} and \mbox{$\Delta^i \psi^j$},
\mbox{$i=2,\dots,L-1$}, by carefully examining
\mbox{$\psi + \pt_t u_h$}, the last term in~\eqref{error-rep-discr}.
For \mbox{$L=2$}, we have
\begin{gather}\label{psi-ptuh-2}
  \psi+\pt_t u_h
    = \Delta^0 \psi^j + \frac{2\Delta^1 u_h^j}{\tau_j}
         + \Delta^1\psi^j P_1\circ \xi_j
       \ \ \ \text{in} \ I_j\,, \ \ j=1,\dots,M,
\end{gather}
while for \mbox{$L>2$}
\begin{gather*}
  \begin{split}
  \psi+\pt_t u_h
    &  = \Delta^0 \psi^j + \frac{2\Delta^1 u_h^j}{\tau_j}
            - \frac{\Delta^2 \psi^j}{3}
         + \sum_{i=1}^{L-3}
           \left\{ \frac{\Delta^i\psi^j}{2i-1}
                     + \frac{2}{\tau_j} \Delta^{i+1} u_h^j
                     - \frac{\Delta^{i+2} \psi^j}{2i+3}
           \right\} P_i\circ \xi_j \\
    & \qquad
         + \left\{ \frac{\Delta^{L-2}\psi^j}{2L-5}
                     + \frac{2}{\tau_j} \Delta^{L-1} u_h^j
           \right\} P_{L-2}\circ \xi_j
         + \frac{\Delta^{L-1}\psi^j}{2L-3} P_{L-1}\circ \xi_j
       \ \ \ \text{in} \ I_j\,, \ \ j=1,\dots,M.
  \end{split}
\end{gather*}
This gives rise to defining the $\Delta^i\psi^j$ as follows
\begin{gather*}
  \frac{\Delta^2 \psi^j}{3} \coloneqq \Delta^0 \psi^j + \frac{2\Delta^1 u_h^j}{\tau_j}
  \quad\text{and}\quad
     \frac{\Delta^{i+2} \psi^j}{2i+3} \coloneqq
     \frac{\Delta^i\psi^j}{2i-1} + \frac{2}{\tau_j} \Delta^{i+1} u_h^j\,, \ \ 
         \ i=1, \dots, L-3, \ j=1,\dots,M.
\end{gather*}
The $\Delta^i u_h^j$ are then computed according to~\eqref{reconstr-delta-psi}.
With these definitions of the quantities $\Delta^i\psi^j$ and $\Delta^iu_h^j$,
we have
\begin{gather}\label{psi-ptuh-L}
  \begin{split}
  \psi+\pt_t u_h
         = \left\{ \frac{\Delta^{L-2}\psi^j}{2L-5}
                     + \frac{2}{\tau_j} \Delta^{L-1} u_h^j
           \right\} P_{L-2}\circ \xi_j
         + \frac{\Delta^{L-1}\psi^j}{2L-3} P_{L-1}\circ \xi_j
       \ \ \ \text{in} \ I_j\,, \ \ j=1,\dots,M.
  \end{split}
\end{gather}

Both~\eqref{psi-ptuh-2} and~\eqref{psi-ptuh-L} take the form
\begin{gather*}
  \begin{split}
  \psi+\pt_t u_h
         = \alpha_L P_{L-2}\circ \xi_j
         + \beta_L P_{L-1}\circ \xi_j \ \ \ \text{in} \ I_j\,, \ \ j=1,\dots,M,
  \end{split}
\end{gather*}
where $\alpha_L$ and $\beta_L$ may vary in space, but are constant in time.
The following lemma will be used to bound those terms.

\begin{lemma}
  Let $v\in V$ be arbitrary.
  Then
  \begin{gather*}
    \abs{\int_{I_j} \dual{v}{\m{G}(T-s)}P_i\bigl(\xi_j(s)\bigr) \D s}
       \le \sigma_j \mu_{j,i} \norm{v}_{\infty,\Omega}\,, \ \ \ i=0,1,\dots, \ \ j=1,\dots,M,
  \end{gather*}
  where the $\sigma_j$ have been defined in~\eqref{def:sigma_j} and
  \begin{gather*}
    \mu_{j,i} = \begin{cases}
                   \kappa_0 \tau_j\,, & \text{for} \ \ i=0, \\[1ex]
                   \displaystyle
                   \min
                     \left\{ \frac{\kappa_0 \tau_j}{2} \norm{P_i}_{1,(-1,1)} \,,
                             \frac{\kappa_1}{\tau_j} \int_{I_j} \frac{\bigl(t_j-t\bigr)\bigl(t-t_{j-1}\bigr)}{T-t} \D t
                        + \kappa_1' \frac{\tau_j^2}{24}
                   \right\} ,
                     & \text{for} \ \ i=1,2,\dots
                \end{cases}
  \end{gather*}
\end{lemma}
\begin{proof}
  First, the H\"older inequality and~\eqref{source:ass} yield
  \begin{gather}\label{intvGP-1}
    \abs{\int_{I_j} \dual{v}{\m{G}(T-s)}P_i\bigl(\xi_j(s)\bigr) \D s}
       \le \norm{v}_{\infty,\Omega} \int_{I_j} \phi_0(T-s) \abs{P_i\bigl(\xi_j(s)\bigr)} \D s
       \le \norm{v}_{\infty,\Omega} \frac{\kappa_0 \sigma_j \tau_j}{2} \norm{P_i}_{1,(-1,1)}\,.
  \end{gather}
  Note that $P_0\equiv 1$ and therefore, $\norm{P_0}_{1,(-1,1)}=2$.

  If $i\ge 1$, then integration by parts gives
  \begin{align*}
    & \int_{I_j} \dual{v}{\m{G}(T-s)}P_i\bigl(\xi_j(s)\bigr) \D s
      = \frac{\tau_j}{2} \int_{I_j} \dual{v}{\m{G}'(T-s)}N_{i+1}\bigl(\xi_j(s)\bigr) \D s \\
    & \qquad
      = \frac{2}{\tau_j} \frac{1}{i(i+1)}
           \int_{I_j} \dual{v}{\m{G}'(T-t)}P_i'\bigl(\xi_j(t)\bigr) \bigl(t-t_j\bigr)\bigl(t-t_{j-1}\bigr) \D t\,.
  \end{align*}
  Thus --- by H\"older's inequality and~\eqref{source:ass} ---
  \begin{align*}
    \abs{\int_{I_j} \dual{v}{\m{G}(T-s)}P_i\bigl(\xi_j(s)\bigr) \D s}
    & \le \frac{2}{\tau_j} \norm{v}_{\infty,\Omega}
                \frac{\norm{P_i'\circ \xi_j}_{\infty,I_j}}{i(i+1)}
           \int_{I_j} \phi_1(T-s) \bigl(t_j-s\bigr)\bigl(s-t_{j-1}\bigr) \D s \\
    & = \frac{1}{\tau_j} \norm{v}_{\infty,\Omega}
           \int_{I_j} \phi_1(T-s) \bigl(t_j-s\bigr)\bigl(s-t_{j-1}\bigr) \D s\,,
  \end{align*}
  because $\norm{P_i'}_{\infty,(-1,1)} = i(i+1)/2$.
  Next, by~\eqref{source:ass}
  \begin{gather}\label{intvGP-2}
    \abs{\int_{I_j} \dual{v}{\m{G}(T-s)}P_i\bigl(\xi_j(s)\bigr) \D s}
      \le \norm{v}_{\infty,\Omega}
            \sigma_j \left(  \frac{\kappa_1}{\tau_j} \int_{I_j}
                             \frac{\bigl(t_j-s\bigr)\bigl(s-t_{j-1}\bigr)}{T-s} \D s
                        + \kappa_1' \frac{\tau_j^2}{24}
                     \right)
  \end{gather}
  Taking the minimum of~\eqref{intvGP-1} and~\eqref{intvGP-2}, we obtain the
  desired result.
\end{proof}

\begin{remark}
  Note that
  \begin{gather*}
    \int_{I_j} \frac{(t_j-t)(t-t_{j-1})}{T-t} \D t
      = \begin{cases}
          \displaystyle
          \frac{\tau_j^2}{2} + \left(T-t_j\right)
                \left[\tau_j +\left(T-t_{j-1}\right)\ln\left(1+\frac{\tau_j}{T-t_j}\right)\right]\,,
          & \ \ j=1,\dots,M-1, \\[1ex]
          \displaystyle
          \frac{\tau_M^2}{2}\,, & \ \ j=M.
        \end{cases}
  \end{gather*}
\end{remark}

%Because of linearity, we have
%%
%\begin{gather}\label{reconstr-delta-esti}
%  \norm{\delta_t E^j}_\infty =
%  \norm{\delta_t\left(u_h - R\right)^j}_\infty
%    \le \eta_\mathrm{ell,\delta}^j
%    \coloneqq \eta\left(\delta_t u_h^j, \delta_t\left(f+\psi\right)^j\right)\,,
%      \ \ j=1,\dots,M.
%\end{gather}

%{\color{blue}

\subsubsection{Error contributions by the spatial discretisation}
\label{sssect:apost-para-spatial}

One term in~\eqref{error-rep-discr} still needs our attention:
\begin{gather*}
  \sum_{j=1}^M \int_{I_j} \dual{\mL E(s)}{\m{G}(T-s)} \D s\,.
\end{gather*}
It captures the error contributions arising from the spatial discretisation.

The properties of the Green's function yield
\begin{gather*}
  \dual{\mL E(s)}{\m{G}(T-s)} = \dual{E(s)}{\mL^* \m{G}(T-s)}
                              = - \dual{E(s)}{\pt_t \m{G}(T-s)}\,,
\end{gather*}
which we integrate by parts to give
\begin{gather*}
  \sum_{j=1}^M \int_{I_j} \dual{\mL E(s)}{\m{G}(T-s)} \D s
       = \dual{E^M}{\m{G}(0)} - \dual{E^0}{\m{G}(T)}
            - \sum_{j=1}^M \int_{I_j} \dual{\pt_t E(s)}{\m{G}(T-s)} \D s\,.
\end{gather*}
H\"older's inequality and \eqref{source:ass} give
\begin{gather}
  \label{est-ell0M}
  \abs{\dual{E^0}{\m{G}(T)}} + \abs{\dual{E^M}{\m{G}(0)}}
     \le \kappa_0 \left\{\sigma_0 \eta_\mathrm{ell}^0 +  \eta_\mathrm{ell}^M
                  \right\}
\end{gather}
Next, note that
\begin{gather*}
  \pt_t E = \frac{2}{\tau} \sum_{i=1}^{L-1} \Delta^i E^j P_{i-1}\circ\xi_j
    \ \ \ \text{in} \ I_j, \ j=1,\dots,M.
\end{gather*}
Consequently, using~\eqref{source:ass} and~\eqref{reconstr-esti-delta}
\begin{gather*}
  \abs{\int_{I_j} \dual{\pt_t E(s)}{\m{G}(T-s)} \D s}
     \le \frac{2}{\tau_j}
         \sum_{i=1}^{L-1} \norm{\Delta^i E^j}_{\infty,\Omega}
                          \int_{I_j} \phi_0(T-s)\abs{P_{i-1}\circ\xi_j} \D s
     \le \kappa_0 \sigma_j
         \sum_{i=1}^{L-1} \eta^j_\mathrm{ell,i} \norm{P_{i-1}}_{1,(-1,1)}\,.
\end{gather*}

%\footnote{
%\begin{gather*}
%  \abs{\dual{\left(\psi+\pt_t u_h\right)(t)}{\m{G}(T-t)} \D t}
%    \le \mu_{j,p-1} \norm{\frac{\Delta^{p-1}\psi^j}{2p-3}
%                     + \frac{2}{\tau_j} \Delta^{p} u_h^j}_{\infty,\Omega}
%        + \mu_{j,p} \norm{\frac{\Delta^{p}\psi^j}{2p-1}}_{\infty,\Omega}
%\end{gather*}
%}
%\footnote{
%and
%%
%\begin{gather*}
%  \abs{\dual{\left(\psi+\pt_t u_h\right)(t)}{\m{G}(T-t)} \D t}
%    \le \sigma_j \tau_j \norm{\Delta^0 \psi^j + \frac{2\Delta^1 u_h^j}{\tau_j}}_{\infty,\Omega}
%         + \mu_{j,1} \norm{\Delta^1\psi^j}_{\infty,\Omega}\,, \ \ j=1,\dots,M,
%\end{gather*}
%
%$R-u_h$\dots{}  see earlier paper.
%}

\subsubsection{A posteriori error bound}
\label{sssect:apost-para-main}

We gather the results of \autoref{sssect:apost-para-gen},
\autoref{sssect:apost-para-temp_rec} and \autoref{sssect:apost-para-spatial}
to establish our main result.

\begin{theorem}\label{theo:theo}
  Let \mbox{$u^j_h\in V_h$}, \mbox{$j=0,1,\dots,M$}, be the sequence of
  approximations to $u(t_j)$ generated by the $L$-step extrapolation
  method introduced in \autoref{sect:discr}.
  Then
  \begin{gather*}
    \norm{u(T)-u_h^M}_{\infty,\Omega}
      \le \eta^M
      \coloneqq \eta_{\mathrm{init}} + \eta_f
        + \eta_{\psi} + \eta_{\Psi},
        + \eta_\mathrm{ell}
  \end{gather*}
  where the components of the error estimator $\eta^M$ are
  \begin{gather*}
    \eta_{\mathrm{init}} \coloneqq
      \kappa_0 \sigma_0 \norm{u^0-u_h^0}_{\infty,\Omega}\,, \quad
    \eta_f \coloneqq
      \kappa_0 \sum_{j=1}^m \sigma_j
                            \int_{I_j} \norm{\bigl(f-\hat{f}\bigr)(s)}_{\infty,\Omega} \D s\,, \\
    \eta_t \coloneqq
      \sum_{j=1}^M \sigma_j \left( \mu_{j,L-2} \norm{\frac{\Delta^{L-2}\psi^j}{2L-5}
                     + \frac{2}{\tau_j} \Delta^{L-1} u_h^j}_{\infty,\Omega}
                     + \frac{\mu_{j,L-1}}{2L-3} \norm{\Delta^{L-1}\psi^j}_{\infty,\Omega} \right) \\
    \intertext{and}
    \eta_\mathrm{ell} \coloneqq
        \kappa_0 \left\{\sigma_0 \eta_\mathrm{ell}^0 +  \eta_\mathrm{ell}^M \right\}
        + \kappa_0 \sum_{j=1}^M \sigma_j
                                \sum_{i=1}^{L-1} \eta^j_\mathrm{ell,i} \norm{P_{i-1}}_{1,(-1,1)}\,.
  \end{gather*}
\end{theorem}

\begin{remark}\label{rem:est}
  \emph{(i)}
  In general, the supremum norms in these estimates can not be determined exactly,
  but need to be approximated.
  For example, one can evaluate $\abs{u^0-u_h^0}$ at a fixed number of points in
  each subdomain of the triangulation of $\Omega$
  (excluding interpolation of quadrature points).

  \emph{(ii)}
  The integrals $\int_{I_j} \norm{\bigl(f-\hat{f}\bigr)(s)}_{\infty,\Omega} \D s$
  need to be approximated as well.
  One possibility is the use of sufficiently accurate quadrature rules, that
  take into account the polynomial degree of the interpolant $\hat{f}$.
\end{remark}

%Direct corollary of this theorem is a posteriori error estimation for semidicertizaction in the time. 
%%
%\begin{corollary}
%  Let $U^M$ is an approximation obtained by extrapolated Euler method in the time $T$. Then, next bound is satisfied
%  %
%  \begin{gather*}
%    \norm{u(T)-U^M}_{\infty,\Omega}
%      \le  \sum_{j=1}^{M} \E^{-\gamma (T-t_j)}
%                  \left(\kappa_0 \eta_F^j
%         +\eta_{\delta\psi_u}^j+\eta_{z_h}^j\right).
%  \end{gather*}
%  %
%\end{corollary}

\bibliographystyle{plain}
% \bibliography{biblio.bib}

\def\cprime{$'$}

\end{document}